\documentclass[10pt,a4paper,oneside,reqno]{amsart}
\usepackage[cp1250]{inputenc}

\title{Special topic}                   
\date{\today}
\usepackage{amssymb}
\usepackage{amsthm} 
\usepackage{esint} 
\usepackage{amsmath}
\usepackage{amsfonts}
\usepackage{mathtools}
\usepackage{amsthm}
\usepackage{amscd}
\usepackage{float}
\usepackage[usenames]{color}
\usepackage[T1]{fontenc}
\usepackage[dvips]{graphicx}
\usepackage{caption} 
\usepackage{algpseudocode}
\usepackage{setspace}

\theoremstyle{definition}
\newtheorem{df}{Definition}
\newtheorem{ex}{Example}
\newtheorem{lem}{Lemma}

\newtheorem{theorem}{Theorem}

\newcommand{\RR}{\mathbb{R}}

\newcommand{\eqnb}{\begin{equation}}
\newcommand{\eqnbs}{\begin{equation*}}
\newcommand{\eqnbsa}{\begin{equation*}\begin{aligned}}
\newcommand{\eqnba}{\begin{equation}\begin{aligned}}
\newcommand{\eqnbl}[1]{\begin{equation}\label{#1}}
\newcommand{\eqnbal}[1]{\begin{equation}\label{#1}\begin{aligned}}
\newcommand{\eqnes}{\end{equation*}}
\newcommand{\eqne}{\end{equation}}

\title{The Lagrange multiplier and the stationary Stokes equations}
\date{21 Mar 2017}

\address{Wojciech Ozanski \\
Mathematics Institute, University of Warwick, Coventry CV4 7AL, UK.}
\email{W.S.Ozanski@warwick.ac.uk}
\begin{document}
\setcounter{page}{1}
\maketitle

\centerline{\begin{minipage}{12cm}
\begin{center}
Abstract
\end{center}
We briefly discuss the notion of the Lagrange multiplier for a linear constraint in the Hilbert space setting, and we prove that the pressure $p$ appearing in the stationary Stokes equations is the Lagrange multiplier of the constraint $\mathrm{div}\, u =0$.
\end{minipage}}
\vspace{1cm}
\section*{Introduction}
For the equations modelling incompressible fluid flows it is frequently remarked that the pressure term acts as a Lagrange multiplier enforcing the incompressibility constraint. Here we prove it rigorously in the case of the stationary Stokes equations. Namely, we show that $p$ appearing in the equations is the Lagrange multiplier corresponding to the constraint $\mathrm{div}\, u =0$ in the variational formulation of the equations, see Section \ref{p_in_SE_section}. For this purpose we briefly discuss preliminary concepts and present some simple variational problems which use Lagrange multipliers in the next section. We then generalise the concept of the Lagrange multiplier to the general Hilbert setting (Section \ref{lagrange_multiplier_section}) and apply it to the stationary Stokes equations. 

\section{Preliminaries}\label{prelims_section}
Let $H$ be a Hilbert space, and $J\colon H\to \RR$ a convex functional that is Frechet differentiable (that is for each $u\in H$ there exists $\nabla J(u) \in H^*$ such that for all $v\in H$ $J(v) = J(u) + \langle \nabla J(u) , v-u \rangle + o(||v-u ||)$, where $\langle u^* , u \rangle $ denotes the duality pairing between a linear functional $u^* \in H^*$ and a point $u\in H$, and $o(x)$ is any function such that $o(x)/x \stackrel{x\to 0^+}{\longrightarrow }0$). Let $K$ denote a closed subspace of $H$.
\begin{lem}\label{lemma1}
If $J: H \to \mathbb{R}$ is convex and differentiable at $u\in K$ then
\[
J(u) = \min_{v\in K} J(v)  \, \Leftrightarrow  \, \nabla J(u) \in K^\circ ,
\]
where $K^\circ \coloneqq \{ u^* \in H^* \colon \langle u^* , v \rangle =0 \text{ for } v\in K \}$ denotes the annihilator of $K$.
\begin{proof}$\mbox{}$
\begin{enumerate}
\item[$(\Rightarrow )$] If $u$ is a minimiser of $J$ over $K$ then $\forall v\in K$ s.t. $||v||=1$ and $\forall t > 0$ we have $J(u\pm tv ) \geq J(u)$ and so
\[ t \langle \nabla J (u) , v \rangle  + o(||tv||) \geq 0 \Rightarrow   \langle \nabla J (u) , v \rangle \geq 0 .\] 
Taking $t<0$ instead of $t>0$, we similarly obtain $\langle \nabla J (u) , v \rangle \leq 0$. Hence $\langle \nabla J (u) , v \rangle = 0$ for all $v\in K$, that is $\nabla J(u)\in K^\circ $.
\item[$(\Leftarrow )$] From convexity, we have
\[
J(u+t(v-u)) = J((1-t)u + tv) \leq t J(v) + (1-t) J(u) \qquad \forall v\in K, \, \forall t\in (0,1).
\]
Subtracting $J(u)$ and dividing by $t$, we get
\[ \frac{1}{t}\left( J(u+t(v-u)) - J(u) \right) \leq J(v) - J(u) . \]
The assumption $\nabla J(u) \in K^\circ$ gives $\langle \nabla J (u) , v \rangle =0$ and so the left hand side is equal to $\frac{o(||t(v-u)||)}{t}$. Taking the limit $t\to 0^+$ we get $0\leq J(v) - J(u)$ for all $v\in K$.\qedhere
\end{enumerate}
\end{proof}
\end{lem}
\begin{ex}\label{ex1}
Let $a_1,\ldots , a_M$ be orthonormal vectors in $H$ and let 
\begin{equation}\label{V}
V = \left\lbrace v\in H \, | \, (a_i , v )=0 \, \, \forall i=1,\ldots , M \right\rbrace ,
\end{equation} 
a finite intersection of hyperplanes. Consider a minimisation problem: Find $u\in V$ such that
\[
J(u) = \min_{v\in V} J(v),
\]
where $J:H \to \mathbb{R}$ is convex and differentiable. Then Lemma \ref{lemma1} gives that $u\in V$ is the minimiser if $\nabla J(u) \in V^\perp = \text{span} \{ a_1, \ldots , a_M \} $. Therefore there exist unique $\lambda_i$'s, $i=1,\ldots , M$, such that $\nabla J(u) = \sum_{i=1}^M \lambda_i a_i$. These $\lambda_i$'s are called \emph{Lagrange multipliers}.
\end{ex}
\begin{ex}\label{ex2} (Elliott \cite{optimisation}, p. 87)
Let $A \in \mathbb{R}^{N\times N}$ be a symmetric, positive definite matrix, $b\in \mathbb{R}^N$, and $C\in \mathbb{R}^{M\times N}$, where $M<N$, be of full rank. Consider the minimisation problem
\begin{equation}\label{ex2prob}
\min_{x\in \text{Ker}\, C} J(x),
\end{equation}
where $J(x):= \frac{1}{2} (x,Ax) - (b,x)$. Note this example is a special case of Example \ref{ex1} with $H\coloneqq \mathbb{R}^N$ and $C = [ a_1 , \ldots , a_M ]^T$ and with a special form of $J$. Hence $\nabla J(u) = \sum_{i=1}^M \lambda_i a_i$. However $a_i = C^T e_i$, where $\{ e_i \}_{i=1, \ldots , M}$ is the standard basis of $\mathbb{R}^M$, and a direct computation shows that $\nabla J(x) = Ax -b$. 
Therefore
\eqnb\label{ex2row}
Ax -b = \sum_{i=1}^M \lambda_i C^T e_i = C^T \Lambda ,
\eqne 
where $\Lambda \coloneqq (\lambda_1 , \ldots , \lambda_M )^T$. We call $\Lambda $ the \emph{Lagrange multiplier} of problem \eqref{ex2prob}.

Rewriting the above equality together with the constraint $x\in \text{Ker}\, C$ in a compact form we obtain
\[
\left[  \begin{array}{cc}
A & -C^T \\
C & 0
\end{array}  \right] \left[  \begin{array}{c}
x \\
\Lambda
\end{array}  \right] = \left[  \begin{array}{c}
b \\
0
\end{array}  \right] .
\]
Since $A$ is invertible and $C$ is of a full rank, the solution $x$ to this system exists and is unique. This example illustrates what is the role of the Lagrange multiplier $\Lambda $: it is a ``redundant variable'' which ``fills out the columns of the system'' and hence makes it solvable for $x$.
\end{ex}
\section{The Lagrange multiplier}\label{lagrange_multiplier_section}
We now generalise the examples to a general Hilbert space setting. Let $M$ be another Hilbert space, let $T\colon H\to M^*$ be a bounded linear operator, and let  $T^*\colon M\to H^*$ denote the dual operator of $T$ (that is $\langle T^*q,u \rangle =\langle Tu ,q \rangle $, where $\langle \cdot , \cdot \rangle $ denotes the duality pairing in turn between $H^*$ and $H$ and between $M^*$ and $M$)

\begin{theorem}\label{mainthm}
Suppose that the operator $T$ satisfies the condition
\begin{equation}\label{c}
|| T^* q ||_{H^*} \geq C || q ||_M   \qquad \text{ for all } q\in M
\end{equation}
for some $C>0$ and consider a minimisation problem: Find $u\in \text{Ker}\, T$ such that
\begin{equation}\label{minimise}
J(u) = \min_{v\in \text{Ker} \, T} J(v).
\end{equation}
Then $u\in \text{Ker} \, T$ is a solution to (\ref{minimise}) if and only if there exists $p\in M$ such that
\[ T^*p=\nabla J(u) .\]
Moreover, if such $p$ exists, it is unique.
\end{theorem}
\begin{df}
This $p\in M$ is the \emph{Lagrange multiplier} of the problem (\ref{minimise}).
\end{df}
Note that by the Fundamental Theorem of Mixed Finite Element Method (see, for example, Lemma 4.1 in Girault \& Raviart \cite{gv}, Chapter I) condition \eqref{c} is equivalent to 
\[
\| T v \|_{M^*} \geq C \| v \|_H   \qquad \text{ for all } v\in ( \mathrm{Ker}\, T)^\perp.
\]

Let us also point out that Example \ref{ex1} is a special case of the theorem above by setting $T=P$, where $P\colon H \to V^\perp$ is an orthogonal projection with respect to the inner product of $H$ (here we identify $H^*$ with $H$). The condition \eqref{c} follows for such $T$ by noting that $M=V^\perp = (\mathrm{Ker} \, T)^\perp$ and by writing
\[ \| T^* q \|_H = \sup_{v\in H} \frac{(T^* q, v)}{\| v \|_H} = \sup_{v\in H} \frac{(Pv ,q)}{\| v \|_H} \geq    \frac{(Pq, q)}{\| q \|_H} = \| q\|_H  \qquad \text{ for all } q\in M . \] 
\begin{proof}(of Theorem \ref{mainthm})
\begin{enumerate}
\item[$(\Leftarrow )$] Since
\[
\langle \nabla J(u) , v \rangle = \langle T^* p , v \rangle = \langle T v  , p \rangle = 0 \qquad \text{ for  } v\in \text{Ker}\, T ,
\]
we see, using Lemma \ref{lemma1}, that $u \in \text{Ker}\, T$ is a solution of (\ref{minimise}).
\item[$(\Rightarrow )$] From (\ref{c}) we can see that $T^*$ is injective on its range $\mathcal{R}(T^*)$. Therefore $T^*$ has a bounded inverse $T ^{-*} : \mathcal{R} (T^*) \to M$ and $\| T  ^{-*} \| \leq \frac{1}{C}$. Hence $T^* : M \to \mathcal{R} (T^*) $ is an isomorphism. In particular $\mathcal{R} (T^* )$ is closed in $H^*$. From Banach Closed Range Theorem (see, for example, Yosida \cite{yosida}, pp. 205-208) we get
\[
\mathcal{R} (T^* ) = \left( \text{Ker} \, T \right)^\circ  . \]
Hence, if $u \in \text{Ker}\, T$ is a solution of (\ref{minimise}), then Lemma \ref{lemma1} gives $\nabla J(u) \in (\mathrm{Ker} \, T)^\circ = \mathcal{R} (T^*)$, that is there exists a unique $ p \in M$ such that $\nabla J(u) = T^* p$. \qedhere
\end{enumerate}
\end{proof}
\section{Pressure function in the stationary Stokes equations}\label{p_in_SE_section}
We now turn into the stationary Stokes equations,
\begin{eqnarray*}
-\Delta {u} + \nabla p &=& f \qquad \text{ in } \Omega ,\\
\text{div} \, {u} &=& 0\qquad \text{ in } \Omega ,  \\
u&=&0 \qquad \text{ on } \partial \Omega ,
\end{eqnarray*} 
where $\Omega \subset \RR^n$ is a smooth domain, ${u} : \Omega \to \mathbb{R}^3$ denotes the velocity of the fluid, $p: \Omega \to \mathbb{R}$ denotes the pressure and $f: \Omega \to \mathbb{R}^3$ is the density of forces acting on the fluid (e.g. gravitational force). The steady Stokes equations govern a flow of a steady, viscous, incompressible fluid. The weak formulation of this problem is to find $u\in V$ and $p\in {L}_0^2$ such that
\eqnb\label{stokes}
(\nabla u , \nabla v ) - (p, \mathrm{div} \, v ) = (f,v) \qquad \text{ for all } v\in H^1_0 (\Omega ),
\eqne
where 
\[
V\coloneqq \left\lbrace {v} \in H_0^1 (\Omega ) \colon \mathrm{div} \, {v} =0 \right\rbrace , \quad L^2_0 (\Omega ) \coloneqq \left\lbrace q\in L^2 (\Omega ) \, \colon \, \int_\Omega q  =0 \right\rbrace 
\]
and $(\cdot , \cdot )$ denotes the $L^2$ inner product (for either scalar, vector or matrix functions).
We will show that the problem \eqref{stokes} is equivalent to finding a minimiser $u\in V$ of the problem
\begin{equation}\label{minstokes}
J(\mathbf{u}) = \min_{\mathbf{v} \in V} J(\mathbf{v}),
\end{equation}
where 
\begin{equation}\label{J}
J(\mathbf{v} )\coloneqq \frac{1}{2} (  \nabla {v}, \nabla v ) - (f , {v}) .
\end{equation}
(Note that this formulation does not include $p$.) Moreover we will show that the pressure function $p$ is the Lagrange multiplier of the problem \eqref{minstokes}.

Indeed, letting $H\coloneqq H^1_0(\Omega )$ and $M\coloneqq  L^2_0 (\Omega )$ we see that $J$ is a convex and differentiable functional on $H$ with 
\[\nabla J(v) = (\nabla  {v} ,\nabla ( \cdot ) ) - (f , \cdot ) \in H^*
\]
for all $v$. Furthermore letting
\[
T\colon H \to M^* \cong L^2_0,\qquad \langle Tv , q \rangle \coloneqq (\mathrm{div } \, v , q ) \quad \text{ for } v\in H, q\in M
\]
we see that $T$ is a bounded linear operator and $V=\mathrm{Ker}\, T$. Moreover $T^* \colon M\to H^*$ is such that $\langle T^* q, {v}  \rangle = \langle T {v} , q \rangle = ( \text{div} \, {v} , q ) $ for $q \in M$, ${v} \in H$, that is
\[
T^* q = \nabla q \quad \text{ as an element of } H^* .
\] 
The condition \eqref{c} follows for such $T^*$ from the well-known inequality $\| q \|_{L^2} \leq  C \| \nabla q \|_{H^{-1} }$ for $q\in L^2_0 (\Omega )$ (see, for example Temam \cite{temam}, pp. 10-11). Therefore Theorem \ref{mainthm} gives that $u\in V$ is a solution to the minimisation problem \eqref{minstokes} if and only if there exists $p\in L^2_0 (\Omega )$ such that 
\[ T^* p = \nabla J(u) = (\nabla  {v} ,\nabla ( \cdot ) ) - (f , \cdot ), \]which is simply the weak formulation \eqref{stokes} of the steady Stokes equations holds.

Note also the similarity of the last equality with \eqref{ex2row}.
\section*{Acknowledgement}
The author is supported by EPSRC as part of the MASDOC DTC at the University of Warwick, Grant No. EP/HO23364/1.
\bibliographystyle{plain}
\bibliography{liter}
\end{document}